\documentclass{amsart}
\usepackage{amscd,amsmath}

\newcommand{\ov}{\overline}

\newcommand{\Ker}{\mathop\mathrm{Ker\,}}

\numberwithin{equation}{section}

\newtheorem{pr}{Proposition}[section]
\newtheorem{co}{Corollary}[section]
\newtheorem{lm}{Lemma}[section]

\theoremstyle{definition}
\newtheorem{de}{Definition}[section]

\begin{document}

\title{Biconformal changes of metric and pseudo-harmonic morphisms}
\author{Radu Slobodeanu}
\address{Faculty of Physics, Bucharest University, Romania}
\email{slobir@k.ro}


\maketitle

\begin{abstract}

Pseudo-harmonic morphisms give rise on the domain space to a distribution which 
admits an almost complex structure compatible with the given Riemannian metric. We shall show 
that this property, together with the harmonicity, are preserved by a biconformal change of
the domain metric. The special case of the pseudo-horizontally homothetic harmonic morphisms
is also treated.\\

\end{abstract}

\section{Introduction}
\begin{text}
Pseudo-harmonic morphisms are a special class of harmonic maps into a Hermitian manifold with
the aditional property called
\emph{Pseudo Horizontal Weak Conformality} (PHWC), cf. \cite {cen}, \cite {lub}. This property
states that the following relation holds: 
\begin{equation} \label{PHWC}
[d\varphi \circ d\varphi^{*}, J]=0,
\end{equation}
where $\varphi: (M, g) \longrightarrow (N, J, h)$ and $d\varphi^{*}$ denotes the adjoint map.

The geometric meaning of (1.1) will become more transparent if we remark that the
 differential of a submersion $\varphi$ from a Riemannian manifold $(M, g)$ into a Hermitian
 manifold $(N, J, h)$ induces an almost complex structure on the horizontal bundle,
defined by 
$J_{\mathcal{H}}=d\varphi\vert _{\mathcal{H}}^{-1} \circ J \circ d\varphi \vert_{\mathcal{H}}$.
One can prove that PHWC 
condition is equivalent to the compatibility of $J_{\mathcal{H}}$ with domain metric $g$.

Moreover, if the codomain is K\"ahler, pseudo-harmonic morphisms have a nice description 
similar to harmonic morphisms. According to \cite {lub}, pseudo-harmonic morphisms  
are those maps that pull back local complex-valued holomorphic functions on $N$ to local 
harmonic functions on $M$.

Then, the analogue of horizontal homothety in this context was introduced by M.A. Aprodu,
 M. Aprodu and V. Brinzanescu in \cite {aab}. 
A \emph{Pseudo-Horizontally Homothetic} (PHH) map is a PHWC map $\varphi$ which satisfies:
\begin{equation} \label{PHH}
[d\varphi \circ \nabla^{M}_{X} \circ d\varphi^{*}, J]=0, \forall X\in\Gamma(\mathcal{H})
\end{equation}

In turn, this condition means that $J_{\mathcal{H}}$ is parallel (with respect to 
$\nabla^{\mathcal{H}}$) in horizontal directions, so satisfies a 
\emph{transversal K\"ahler} condition. Any PHH harmonic submersion onto a K\"ahler manifold
 exhibits a particularly nice geometric property: it pulls back complex submanifolds 
into minimal submanifolds, cf. \cite {aab}. Moreover, in \cite {ma} it is shown that 
PHH harmonic
 submersions are (weakly) stable (in particular such maps minimise the energy-functional
in their homotopy class). Further properties and examples of PHH harmonic submersions can be 
found in \cite {aa}, \cite{brz}.

In this paper we shall explore the invariance of the 
above two geometric conditions (PHWC and PHH) together with harmonicity under certain 
changes of metric. These results are particularly interesting in finding new examples of 
pseudo-harmonic morphisms and also they point out a nice analogy with the case of harmonic 
morphisms.

\end{text}

\section{Changes of metric which preserve pseudo-harmonic morphisms}

According to \cite {mo}, for a PHWC map 
$\varphi :(M^{m}, g) \longrightarrow (N^{2n}, J, h)$ from a Riemannian manifold to a 
K\"ahler one, the tension field is given by:
$$
\tau(\varphi)=-d \varphi(F^{\varphi} \delta F^{\varphi}),
$$
where $F^{\varphi}$ is the $f$-structure on $M$, naturally induced by the PHWC map $\varphi$
 (in fact, $F^{\varphi}$ extends $J_{\mathcal{H}}$ with zero on $\mathcal{V}$) and 
$\delta F^{\varphi}$ = trace $\nabla F^{\varphi}$ is the divergence of $F^{\varphi}$.

If we consider an adapted frame $\{e_{i}, F^{\varphi}e_{i}, e_{\alpha}\}$ (i.e. an orthonormal frame such that
$e_{\alpha} \in \Ker F^{\varphi}$), then the local form of the above formula is:
\begin{equation} \label{deftau}
\tau(\varphi)=-d \varphi\left( \sum_{i=1}^{n}F^{\varphi}\left[(\nabla_{e_{i}}F^{\varphi})(e_{i})+
(\nabla_{F^{\varphi}e_{i}}F^{\varphi})(F^{\varphi}e_{i})\right]+(m-2n) \mu^{\mathcal{V}} \right).
\end{equation}

\begin{de}
Let $(M, g)$ be a Riemannian manifold endowed with two complementary distributions 
$\mathcal{V}$
and $\mathcal{H}$ which induce a natural decomposition of the metric: 
$g=g^{\mathcal{V}}+g^{\mathcal{H}}$.
Then one call \emph{biconformal change of metric} $g$, the association of a new metric on 
$M$ defined by:
$$
\ov{g}=\sigma^{-2}g^{\mathcal{H}}+ \rho^{-2}g^{\mathcal{V}},
$$
where $\sigma, \rho :M \longrightarrow (0, \infty)$ are smooth functions.
\end{de}

The next result for PHWC submersions is similar to that one for 
horizontally conformal submersions (see \cite{ud},
Lemma 4.6.6, p. 126).
\begin{pr}
Let $\varphi :(M^{m}, g) \longrightarrow (N^{2n}, J, h)$ be a PHWC submersion. With respect to the vertical 
and horizontal distributions, one take a biconformal change of the metric $g$, denoted by $\ov{g}$.
Then, regarded as a map from $(M, \ov{g})$ to $(N^{2n}, J, h)$, $\varphi$ has tension field:
\begin{equation}\label{tau}
\ov{\tau}(\varphi)=\sigma^{2}\left[\tau(\varphi)+
d\varphi(\emph{grad}(\mathrm{ln} \rho^{2n-m}\sigma^{2-2n}))\right]
\end{equation}
\end{pr}

\begin{proof}
Consider an adapted frame $\{e_{i}, F^{\varphi}e_{i}, e_{\alpha}\}$ for $(M, g)$.
 Then $\{\sigma e_{i}, \sigma F^{\varphi}e_{i}, \rho e_{\alpha}\}$ will be an adapted frame with respect 
to the perturbed metric $\ov{g}$. Let $f_{i}$ denote $e_{i}$ when $i=\ov{1,n}$ and 
$F^{\varphi}e_{i-n}$ when $i=\ov{n+1,2n}$.
Then, using Koszul formula, one obtain for any $X, Y \in\Gamma(\mathcal{H})$:
\begin{equation}\label{Kos1}
\mathcal{H}(\ov{\nabla}_{X}Y)=\mathcal{H}(\nabla_{X}Y)+
\sum_{i=1}^{n}\left[-X(\mathrm{ln} \sigma)g(Y, f_{i})-Y(\mathrm{ln} \sigma)g(X, f_{i})+f_{i}(\mathrm{ln} \sigma)g(X, Y)\right]f_{i}.
\end{equation}

In the same way, but for $V \in\Gamma(\mathcal{V})$ one get:
\begin{equation}\label{Kos2}
\mathcal{H}(\ov{\nabla}_{V}V)=\frac{\sigma^{2}}{2}\left[2\rho^{-2}\mathcal{H}(\nabla_{V}V)+
\sum_{i=1}^{n}\left[-f_{i}(\rho^{-2})g(V, V)\right]f_{i}\right].
\end{equation}

Now, apllying \eqref{Kos2} for $V=\rho e_{\alpha}, \alpha=\ov{2n+1,m}$ and summing, one 
get the following result which has the same form as in the harmonic morphism case 
(see \cite {ud}, Remark 4.6.8, p. 126):
\begin{lm}
At a biconformal change of the metric, the mean curvature of the fibers of a 
pseudo-harmonic morphism changes in this way:
\begin{equation} \label{mu}
\ov{\mu}^{\mathcal{V}}=\sigma^{2}\left[\mu^{\mathcal{V}}+ \mathcal{H}(\mathrm{grad}(\mathrm{ln} \rho))\right]
\end{equation}
\end{lm}

\bigskip
Now, in order to get the formula stated by the theorem, one has to see how changes
 the first term in equation \eqref{deftau}, that is:
\begin{equation*}
\begin{split}
F^{\varphi}\mathrm{div}_{\mathcal{H}}F^{\varphi}=
&F^{\varphi}\left[(\nabla_{e_{i}}F^{\varphi})(e_{i})+
(\nabla_{F^{\varphi}e_{i}}F^{\varphi})(F^{\varphi}e_{i})\right]=\\
&\mathcal{H}(\nabla_{e_{i}} e_{i})+
\mathcal{H}(\nabla_{F^{\varphi} e_{i}}F^{\varphi} e_{i})
+F^{\varphi}(\nabla_{e_{i}}F^{\varphi}e_{i})-F^{\varphi}(\nabla_{F^{\varphi}e_{i}}e_{i}).
\end{split}
\end{equation*}

So we have to apply formula \eqref{Kos1} in four particular cases:

\medskip
$\mathcal{H}(\ov{\nabla}_{\sigma e_{i}}\sigma e_{i})=\sigma^{2}\mathcal{H}(\nabla_{e_{i}}e_{i})
-\sigma^{2}e_{i}(\mathrm{ln} \sigma)e_{i}+
\sigma^{2}\mathrm{grad}_{\mathcal{H}}(\mathrm{ln}\sigma)$,

\medskip
$\mathcal{H}(\ov{\nabla}_{\sigma F^{\varphi}e_{i}}\sigma F^{\varphi}e_{i})=
\sigma^{2}\mathcal{H}(\nabla_{F^{\varphi}e_{i}}F^{\varphi}e_{i})
-\sigma^{2}F^{\varphi}e_{i}(\mathrm{ln} \sigma)F^{\varphi}e_{i}+
\sigma^{2}\mathrm{grad}_{\mathcal{H}}(\mathrm{ln} \sigma)$,

\medskip
$\mathcal{H}(\ov{\nabla}_{\sigma e_{i}}\sigma F^{\varphi}e_{i})=
\sigma^{2}\mathcal{H}(\nabla_{e_{i}}F^{\varphi}e_{i})
-\sigma^{2}F^{\varphi}e_{i}(\mathrm{ln} \sigma)e_{i}$,

\medskip
$\mathcal{H}(\ov{\nabla}_{\sigma F^{\varphi}e_{i}}\sigma e_{i})=
\sigma^{2}\mathcal{H}(\nabla_{F^{\varphi}e_{i}}e_{i})
-\sigma^{2}e_{i}(\mathrm{ln} \sigma)F^{\varphi}e_{i}$.

\medskip

Now, summing these four terms as they appear in $F^{\varphi}\ov{\mathrm{div}}_{\mathcal{H}}F^{\varphi}$ we get:
\begin{equation}\label{div}
F^{\varphi}\ov{\mathrm{div}}_{\mathcal{H}}F^{\varphi}=\sigma^{2}\left(F^{\varphi}\mathrm{div}_{\mathcal{H}}F^{\varphi}+(2n-2)\mathrm{grad}(\mathrm{ln}\sigma)\right)
\end{equation}

Now, from \eqref{mu} and \eqref{div}, the stated result follows.
\end{proof}

The following result is analogous with the one proved in \cite {mox}.
\begin{co}
Let $\varphi :(M^{m}, g) \longrightarrow (N^{2n}, J, h)$ be a PHWC submersion (with $m>2n$).
 For any smooth function $\sigma:M \longrightarrow (0, \infty)$, set
\begin{equation}\label{change}
g_{\sigma}=\sigma^{-2}g^{\mathcal{H}}+ \sigma^{\frac{4n-4}{m-2n}}g^{\mathcal{V}}.
\end{equation}
Then $\varphi$ is a pseudo-harmonic morphism with respect to $g$ if and only if
 is a a pseudo-harmonic morphism with respect to $g_{\sigma}$.
\end{co}

\bigskip
Now, suppose that $\varphi$ is a pseudo-horizontally homothetic map. 
The condition \eqref {PHH} translates simply into 
$d\varphi \left( (\nabla_{X}F^{\varphi})Y \right)=0, \forall X, Y \in \Gamma(\mathcal{H})$, 
so:
$$\left[(\nabla_{X}F^{\varphi})Y\right]^{\mathcal{H}}=0, \forall X, Y \in \Gamma(\mathcal{H}).
$$
In particular, $F^{\varphi}\mathrm{div}_{\mathcal{H}}F^{\varphi}=0$, so the formula
\eqref{deftau} reduces to:
$$
\tau(\varphi)=-(m-2n)d \varphi\left(\mu^{\mathcal{V}} \right).
$$
From the above relation, one can see that a pseudo-horizontally homothetic map is harmonic 
if and only if the fibers are minimal (as it was proved in \cite {aab}).

We want to see under which changes of metric of the type \eqref{change} the PHH harmonic morphisms are invariant.
\begin{co}
Let $\varphi :(M^{m}, g) \longrightarrow (N^{2n}, J, h)$ be a PHH harmonic morphism (with $m>2n$). Then $\varphi$ is a
PHH harmonic morphism with respect to $g_{\sigma}$ if and only if $\sigma$ is constant. In particular,
for any any strictly positive constant $c$, set
$$
g_{c}=c^{-2}g^{\mathcal{H}}+ c^{\frac{4n-4}{m-2n}}g^{\mathcal{V}}.
$$
Then a PHH submersion $\varphi$ is a PHH harmonic morphism with respect to $g$ if and only if
 is a a PHH harmonic morphism with respect to $g_{c}$.
\end{co}

\begin{proof}
In particular, $\varphi$ is a PHWC submersion. Using \eqref{Kos1}, for any $X$ tangent to
$\mathcal{H}$ we derive:
\begin{equation*}
\begin{split}
\mathcal{H}\left((\ov{\nabla}_X F^{\varphi})Y \right)=&
\mathcal{H}\left((\nabla_X F^{\varphi})Y\right)
+\mathrm{grad}_{\mathcal{H}}(\mathrm{ln}\sigma)g(X, F^{\varphi})Y)+\\
&g(X, Y)\sum_{i=1}^{n}\left[F^{\varphi}e_{i}(\mathrm{ln}\sigma)e_{i}+
e_{i}(\mathrm{ln}\sigma)F^{\varphi}e_{i}\right]-\\
&F^{\varphi}Y(\mathrm{ln}\sigma)X+
Y(\mathrm{ln}\sigma)F^{\varphi}X.
\end{split}
\end{equation*}
As $\varphi$ is a PHH harmonic morphism with respect to $g$, if $\sigma$ is constant,
then it is obvious from the above formula that $\varphi$ is PHH harmonic morphism with respect
to $g_{\sigma}$, too.

Conversely, if $\varphi$ is PHH harmonic morphism with respect to $g_{\sigma}$, then
$F^{\varphi}\ov{\mathrm{div}}_{\mathcal{H}}F^{\varphi}=0$. But the hypothesis implies also
$F^{\varphi}\mathrm{div}_{\mathcal{H}}F^{\varphi}=0$ and the relation \eqref{div} tells us 
precisely that $\sigma$ must be a constant. 
\end{proof}

\end{document}